\documentclass[11pt]{article}
\usepackage[russian]{babel}
\usepackage[cp1251]{inputenc}
\usepackage{amsmath,latexsym}
\usepackage[psamsfonts]{amssymb}
\usepackage[dvips]{graphicx}
\usepackage[mathcal]{euscript}
\begin{document}
УДК. 517.954, 517.927.25,517.589
\begin{center}
\textbf{ Нелокальная задача с условиями сопряжения \\ для вырождающегося
уравнения высокого порядка\\  с дробной производной Капуто}\\
\textbf{Б.Ю.Иргашев}\\
Наманганский инженено-строительный институт.г.Наманган.Узбекистан.\\
Институт Математики Республики Узбекистан.\\
bahromirgasev@gmail.com
\end{center}

\textbf{Аннотация.} В работе для вырождающегося уравнения высокого порядка с  дробной производной в смысле Капуто изучена нелокальная задача с условиями сопряжения, в прямоугольной области. Решение построено в виде ряда Фурье по собственным функциям одномерной задачи.  Дан критерий единственности решения.

 \textbf{Ключевые слова.} Уравнение четного порядка, дробная производная Капуто, вырождение, условия сопряжения, собственное значение, собственная функция, ряд Фурье, сходимость.

\textbf{1.Введение и постановка задачи.} В области $\Omega  = {\Omega _x} \times {\Omega _y},$ ${\Omega _x} = \left\{ {x:0 < x < 1 } \right\},$ ${\Omega _y} = \left\{ {y:\,  -a < y < b} \right\}$, $a ,b > 0,\,a,b \in R$, рассмотрим уравнение в частных производных
\[L\left[ u \right] \equiv \left\{ \begin{gathered}
  {}_CD_{0y}^\alpha u\left( {x,y} \right) + {\left( { - 1} \right)^k}{x^m}\frac{{{\partial ^{2k}}u\left( {x,y} \right)}}{{\partial {x^{2k}}}} = 0,\,y > 0,0 < \alpha  < 1, \hfill \\
  {}_CD_{0y}^\beta u\left( {x,y} \right) + {\left( { - 1} \right)^k}{x^m}\frac{{{\partial ^{2k}}u\left( {x,y} \right)}}{{\partial {x^{2k}}}} = 0 = 0,y < 0,1 < \beta  < 2, \hfill \\
\end{gathered}  \right.\eqno(1)\]
где
$$0 \leqslant m < k,\,m \notin N,$$
\[{}_CD_{0y}^\alpha u\left( {x,y} \right) = \frac{1}{{\Gamma \left( {1 - \alpha } \right)}}\int\limits_0^y {\frac{{\frac{{\partial u\left( {x,z} \right)}}{{\partial z}}}}{{{{\left( {y - z} \right)}^\alpha }}}} dz,{}_CD_{0y}^\beta u\left( {x,y} \right) = \frac{1}{{\Gamma \left( {2 - \beta } \right)}}\int\limits_0^y {\frac{{\frac{{{\partial ^2}u\left( {x,z} \right)}}{{\partial {z^2}}}}}{{{{\left| {y - z} \right|}^{\beta  - 1}}}}} dz\]
- дробные производные в смысле Капуто.\\
Пусть ${\Omega _ + } = \Omega  \cap \left( {y > 0} \right),{\mkern 1mu} {\mkern 1mu} {\Omega _ - } = \Omega  \cap \left( {y < 0} \right).$  Для уравнения  (1) рассмотрим следующую  задачу.

\textbf{Задача D.}  Найти функцию $u(x,y)$  с условиями
\[Lu\left( {x,y} \right) \equiv 0,{\mkern 1mu} {\mkern 1mu} {\mkern 1mu} {\mkern 1mu} \left( {x,y} \right) \in {\Omega _ + } \cup {\Omega _ - },\]
\[\frac{{{\partial ^{2k}}u}}{{\partial {x^{2k}}}} \in C\left( \Omega  \right),\frac{{{\partial ^{2k - 1}}u}}{{\partial {x^{2k - 1}}}} \in C\left( {\bar \Omega } \right),{x^{\frac{m}{2}}}\frac{{{\partial ^{2k}}u}}{{\partial {x^{2k}}}} \in {L_2}\left( {{{\overline \Omega  }_x}} \right),y \in \left[ { - a,b} \right],\]
$${}_CD_{0y}^\alpha u \in C\left( {{\Omega _ + }} \right),{}_CD_{0y}^\beta u \in C\left( {{\Omega _ - }} \right),$$
$$\frac{{{\partial ^{j}}u}}{{\partial {x^{j}}}}\left( {0,y} \right) = \frac{{{\partial ^{j}}u}}{{\partial {x^{j}}}}\left( {1 ,y} \right) = 0,{\mkern 1mu} {\kern 1pt} j = 0,1,...,k - 1; - a \leqslant y \leqslant b,\eqno(2)$$
$$u\left( {x, + 0} \right) = u\left( {x, - 0} \right),\,  _CD_{0y}^\alpha u\left( {x, + 0} \right) = {u_y}\left( {x, - 0} \right),0 \leqslant x \leqslant 1,\eqno(3)$$
\[u\left( {x,b} \right) - u\left( {x,-a} \right) = \varphi \left( x \right),\,0 \leqslant x \leqslant 1,\eqno(4)\]
где $\varphi \left( x \right)$ достаточно гладкая функция.

Дифференциальные уравнения в частных производных дробного порядка лежат в основе математического моделирования различных физических процессов и явлений окружающей среды, имеющих фрактальную природу [1], [2].

Уравнение (1) при $\alpha =1$, $\beta =2$, $k =1$ является параболо-гиперболическим и исследовалось во многих работах, например [3]-[9].
На необходимость рассмотрения задач сопряжения, когда на одной части области задано параболическое уравнение, на другой – гиперболическое, было указано в 1959 году И.М.Гельфандом [3]. В этой работе приводится пример, связанный с движением газа в канале, окруженном пористой средой: в канале движение газа описывается волновым уравнением, вне его – уравнением диффузии. Одной из первых работ, посвященных изучению краевых задач для параболо-гиперболических уравнений, явилась работа Г.М.Стручиной [4]. Затем Я.С.Уфлянд [5] задачу о распространения электрических колебаний в составных линиях, когда на участке   полу бесконечной линии пренебрегается потерями, а остальная часть линии рассматривается как кабель без утечки, свел к решению краевых задач для параболо-гиперболических уравнений.

Задачи с условиями сопряжения для уравнений с дробными производными изучались в работах [10]-[14] для уравнений второго порядка, а в работах [15],[16] для уравнений четвертого порядка c нелокальными условиями.

В данной работе изучается краевая задача с нелокальным условием и условиями сопряжения для уравнения высокого четного порядка c вырождением. Решение построено в виде ряда Фурье. Получены достаточные условия возможности почленного дифференцирования ряда. Найдены условия единственности решения.

\textbf{2. Существование решения.} Заметим, что из условия (2) и применения теоремы Лагранжа следует, что
\[u\left( {x,y} \right) = O\left( {{x^k}} \right),\,x \to  + 0.\]
 Будем искать решение поставленной задачи методом Фурье
\[u\left( {x,y} \right) = X\left( x \right)Y\left( y \right).\eqno(5)\]
Подставим (5) в уравнение (1). Тогда переменные разделятся . Учитывая условия (2), относительно переменной $x$ , получим задачу на собственные значения
\[\left\{ \begin{gathered}
  {X^{\left( {2k} \right)}}\left( x \right) = {\left( { - 1} \right)^k}\lambda {x^{ - m}}X\left( x \right), \hfill \\
  {X^{\left( j \right)}}\left( 0 \right) = {X^{\left( j \right)}}\left( 1 \right) = 0,\,j = 0,1,...,k - 1. \hfill \\
\end{gathered}  \right.\eqno(6)\]
Ясно,что $\lambda  = 0$ не является собственным значением. Покажем, что  $\lambda  > 0.$ Действительно, имеем
\[\int\limits_0^1 {X\left( x \right){X^{\left( {2k} \right)}}\left( x \right)dx}  = {\left( { - 1} \right)^k}\lambda \int\limits_0^1 {{x^{ - m}}{X^2}\left( x \right)dx} ,\]
используя формулу интегрирования по частям несколько раз, получим
\[\int\limits_0^1 {{{\left\{ {{X^{\left( k \right)}}\left( x \right)} \right\}}^2}dx}  = \lambda \int\limits_0^1 {{x^{ - m}}{X^2}\left( x \right)dx} ,\]
отсюда
$$\lambda  > 0.$$
Покажем, что задача (6) имеет собственные значения . Для этого приведем эту задачу к эквивалентному интегральному  уравнению с помощью функции Грина. Интегральное уравнение имеет вид
\[X\left( x \right) = {\left( { - 1} \right)^k}\lambda \int\limits_0^1 {{\xi ^{ - m}}G\left( {x,\xi } \right)X\left( \xi  \right)d\xi },\eqno(7)\]
где
\[G\left( {x,\xi } \right) =  - \frac{1}{{\left( {2k - 1} \right)!}}\left\{ \begin{gathered}
  {G_1}\left( {x,\xi } \right),\,0 \leqslant x \leqslant \xi , \hfill \\
  {G_2}\left( {x,\xi } \right),\,\xi  \leqslant x \leqslant 1, \hfill \\
\end{gathered}  \right.\]
- функция Грина задачи (6) (см. [17]), здесь
$${G_1}\left( {x,\xi } \right) = {\left( {1 - \xi } \right)^k}{x^k}\sum\limits_{i = 0}^{k - 1} {\sum\limits_{j = 0}^{k - i - 1} {{{\left( { - 1} \right)}^i}C_{2k - 1}^i} } C_{k - 1 + j}^j{x^{k - i - 1}}{\xi ^{j + i}},$$
$${G_2}\left( {x,\xi } \right) = {\left( {1 - x} \right)^k}{\xi ^k}\sum\limits_{i = 0}^{k - 1} {\sum\limits_{j = 0}^{k - i - 1} {{{\left( { - 1} \right)}^i}C_{2k - 1}^i} } C_{k - 1 + j}^j{\xi ^{k - i - 1}}{x^{j + i}},$$
$$C_n^k = \frac{{n!}}{{k!\left( {n - k} \right)!}}.$$
Запишем (7) в виде
$${x^{ - \frac{m}{2}}}X\left( x \right) = \lambda \int\limits_0^1 {{\xi ^{ - \frac{m}{2}}}\left[ {{{\left( { - 1} \right)}^k}G\left( {x,\xi } \right)} \right]{x^{ - \frac{m}{2}}}\left( {{x^{ - \frac{m}{2}}}X\left( \xi  \right)} \right)d\xi }.$$
Введем обозначения
\[\overline X \left( x \right) = {x^{ - \frac{m}{2}}}X\left( x \right),\]
$$\overline G \left( {x,\xi } \right) = {\xi ^{ - \frac{m}{2}}}\left[ {{{\left( { - 1} \right)}^k}G\left( {x,\xi } \right)} \right]{x^{ - \frac{m}{2}}},$$
тогда имеем уравнение
$$\overline X \left( x \right) = \lambda \int\limits_0^1 {\overline G \left( {x,\xi } \right)\overline X \left( x \right)dx}.\eqno(8)$$
(8) - есть интегральное уравнение с непрерывным, по обоим переменным, и симметричным ядром. По теории уравнений с симметричными ядрами уравнение (8) имеет не более чем счетное число собственных значений и собственных функций. Итак, задача (6) имеет собственные значения
$0 < {\lambda _1} < {\lambda _2} < ....,\,\mathop {\lim }\limits_{n \to  + \infty } {\lambda _n} = \infty ,$ и соответствующие им собственные функции -$ {X_1}(x),{X_2}\left( x \right),... .$ Далее будем считать, что
\[{\left\| {{X_n}(x)} \right\|^2} = \int\limits_0^1 {{x^{ - m}}X_n^2(x)dx}  = 1,\]
тогда справедливо неравенство Бесселя
$$\sum\limits_{n = 1}^\infty  {{{\left( {\frac{{{X_n}\left( x \right)}}{{{\lambda _n}}}} \right)}^2}}  \leqslant \int\limits_0^1 {{\xi ^{ - m}}{G^2}\left( {x,\xi } \right)d\xi }  < \infty .\eqno(9)$$
Теперь найдем условия, при которых заданная функция $\varphi \left( x \right)$  разлагается в ряд по собственным функциям ${X_n}(x)$. Для этого воспользуемся теоремой Гильберта-Шмидта.

\textbf{Теорема 1}. Пусть функция $\varphi \left( x \right)$ удовлетворяет условиям:\\
1. $ {x^{ - \frac{m}{2}}}\varphi \left( x \right) \in C\left[ {0;1} \right];$\\
2. $\varphi \left( x \right) \in {C^{2k}}\left[ {0,1} \right];$\\
3. ${\varphi ^{\left( i \right)}}\left( 0 \right) = {\varphi ^{\left( i \right)}}\left( 1 \right) = 0,\,i = 0,1,...,k - 1.$\\
Тогда её  можно разложить в равномерно сходящий ряд вида
\[\varphi \left( y \right) = \sum\limits_{n = 1}^\infty  {{\varphi _n}{X_n}\left( x \right)} ,\]
где
\[{\varphi _n} = \int\limits_0^1 {{x^{ - m}}\varphi \left( x \right){X_n}(x)dx} .\]

\textbf{Доказательство.} Покажем справедливость равенства
\[{x^{ - \frac{m}{2}}}\varphi \left( x \right) = \int\limits_0^1 {\overline G \left( {x,\xi } \right)\left( {{{\left( { - 1} \right)}^k}{\xi ^{\frac{m}{2}}}\frac{{{d^{2k}}\varphi \left( \xi  \right)}}{{d{\xi ^{2k}}}}} \right)d\xi }.\]
Действительно
\[\int\limits_0^1 {{\xi ^{ - \frac{m}{2}}}\left[ {{{\left( { - 1} \right)}^k}G\left( {x,\xi } \right)} \right]{x^{ - \frac{m}{2}}}\left( {{{\left( { - 1} \right)}^k}{\xi ^{\frac{m}{2}}}\frac{{{d^{2k}}\varphi \left( \xi  \right)}}{{d{\xi ^{2k}}}}} \right)d\xi }  = \]
\[ = {x^{ - \frac{m}{2}}}\int\limits_0^1 {G\left( {x,\xi } \right)} \frac{{{d^{2k}}\varphi \left( \xi  \right)}}{{d{\xi ^{2k}}}}d\xi  = {x^{ - \frac{m}{2}}}\varphi \left( x \right).\]
Т.е. для функции ${x^{ - \frac{m}{2}}}\varphi \left( x \right)$  выполняются условия теоремы Гильберта-Шмидта и значит
\[{x^{ - \frac{m}{2}}}\varphi \left( x \right) = \sum\limits_{n = 1}^\infty  {{x^{ - \frac{m}{2}}}{\varphi _n}{X_n}\left( x \right)} ,\]
разделив  на ${x^{ - \frac{m}{2}}}$ , имеем
\[\varphi \left( x \right) = \sum\limits_{n = 1}^\infty  {{\varphi _n}{X_n}\left( x \right)} .\]
\textbf{Теорема 1} доказана.

Перейдем к решению краевой задачи для уравнения по переменной $y$. Учитывая (5), (1) и условия (3),(4) получим
\[\left\{ \begin{gathered}
  {}_CD_{0y}^\alpha {Y_n}\left( y \right) =  - {\lambda _n}{Y_n}\left( y \right),\,y > 0, \hfill \\
  {}_CD_{0y}^\beta {Y_n}\left( y \right) =  - {\lambda _n}{Y_n}\left( y \right),y < 0, \hfill \\
  {Y_n}\left( { + 0} \right) = {Y_n}\left( { - 0} \right), \hfill \\
  {}_CD_{0y}^\alpha {Y_n}\left( { + 0} \right) = {{Y'}_n}\left( { - 0} \right), \hfill \\
  {Y_n}\left( b \right) - {Y_n}\left( -a \right) = {\varphi _n} = \int\limits_0^1 {\varphi \left( x \right){X_n}\left( x \right)dx} . \hfill \\
\end{gathered}  \right.\eqno(10)\]
Решение уравнений входящих в задачу (10) будем искать в виде
\[{Y_n}\left( y \right) = \left\{ \begin{gathered}
  \sum\limits_{s = 0}^\infty  {{c_{1s}}{y^{{\gamma _1}s + {\delta _1}}}} ,\,y > 0, \hfill \\
  \sum\limits_{s = 0}^\infty  {{c_{2s}}{{\left( { - y} \right)}^{{\gamma _2}s + {\delta _2}}}} ,\,y < 0, \hfill \\
\end{gathered}  \right.\eqno(11)\]
где ${c_{1s}},{c_{2s}},{\gamma _1},{\gamma _2},{\delta _1},{\delta _2}$ неизвестные постоянные подлежащие определению. Подставляя (11) в уравнение (10) и приравнивая коэффициенты при одинаковых степенях $y$, находим общее решение уравнения (10) в виде
\[{Y_n}\left( y \right) = {c_1}{E_{\alpha ,1}}\left( { - {\lambda _n}{y^\alpha }} \right),\,(y > 0),\]
\[{Y_n}\left( y \right) = {c_2}{E_{\beta ,1}}\left( { - {\lambda _n}{{\left( { - y} \right)}^\beta }} \right) + {c_3}\left( { - y} \right){E_{\beta ,2}}\left( { - {\lambda _n}{{\left( { - y} \right)}^\beta }} \right),\,(y < 0),\]
где
\[{c_1},{c_2},{c_3} - const,\]
\[{E_{\mu ,\eta }}\left( z \right) = \sum\limits_{n = 0}^\infty  {\frac{{{z^n}}}{{\Gamma \left( {\mu n + \eta } \right)}}}\eqno(12)\]
-	функция Миттаг-Леффлера.\\
Удовлетворив условиям задачи (10) имеем
\[{c_1} = {c_2},\,{\lambda _n}{c_1} = {c_3},\]
\[{c_1} = \frac{{{\varphi _n}}}{{\Delta \left( n \right)}},\]
где
\[\Delta \left( n \right) = {E_{\alpha ,1}}\left( { - {\lambda _n}{b^\alpha }} \right) - \left( {{E_{\beta ,1}}\left( { - {\lambda _n}{a^\beta }} \right) + a{\lambda _n}{E_{\beta ,2}}\left( { - {\lambda _n}{a^\beta }} \right)} \right).\]
Отсюда решение задачи (10) имеет вид
\[{Y_n}\left( y \right) = \frac{{{\varphi _n}}}{{\Delta \left( n \right)}}{E_{\alpha ,1}}\left( { - {\lambda _n}{y^\alpha }} \right),\,(y > 0),\]
\[{Y_n}\left( y \right) = \frac{{{\varphi _n}}}{{\Delta \left( n \right)}}{E_{\beta ,1}}\left( { - {\lambda _n}{{\left( { - y} \right)}^\beta }} \right) + \frac{{{\lambda _n}{\varphi _n}}}{{\Delta \left( n \right)}}\left( { - y} \right){E_{\beta ,2}}\left( { - {\lambda _n}{{\left( { - y} \right)}^\beta }} \right),\,(y < 0).\]

Теперь нам надо показать отделимость от нуля выражения $\Delta \left( n \right)$ . Известно [18] , что функция (11) на вещественной оси, если
\[0 < \mu  < 2,\]
может иметь только конечное число нулей. В частности, если
\[1 < \mu  < 2,\,\eta  \geqslant \frac{3}{2}\mu ,\]
функция (11) не имеет вещественных нулей [19].\\
Справедлива следующая лемма.

\textbf{Лемма .} Пусть
$$h = \max \left\{ {t:\,{E_{\beta ,2}}\left( { - t} \right) = 0,t > 0} \right\},$$
тогда существует номер $N > {n_0},$ где ${\lambda _{{n_0}}}{a^\beta } > h,$ что для всех $n > N$
выполняется оценка
\[\left| {{E_{\alpha ,1}}\left( { - {\lambda _n}{b^\alpha }} \right) - \left( {{E_{\beta ,1}}\left( { - {\lambda _n}{a^\beta }} \right) + a{\lambda _n}{E_{\beta ,2}}\left( { - {\lambda _n}{a^\beta }} \right)} \right)} \right| \geqslant \delta  > 0,\]
где $\delta $  некоторое положительное число.

\textbf{Доказательство. } Используя асимптотические разложения из [18] , при ${\lambda _n} \to  + \infty $,  имеем
\[{E_{\alpha ,1}}\left( { - {\lambda _n}{b^\alpha }} \right) = \frac{1}{{{\lambda _n}{b^\alpha }\Gamma \left( {1 - \alpha } \right)}} + O\left( {\frac{1}{{\lambda _n^2}}} \right),\]
\[{E_{\beta ,1}}\left( { - {\lambda _n}{a^\beta }} \right) = \frac{1}{{{\lambda _n}{a^\beta }\Gamma \left( {1 - \beta } \right)}} + O\left( {\frac{1}{{\lambda _n^2}}} \right),\]
\[{E_{\beta ,2}}\left( { - {\lambda _n}{a^\beta }} \right) = \frac{1}{{{\lambda _n}{a^\beta }\Gamma \left( {2 - \beta } \right)}} + O\left( {\frac{1}{{\lambda _n^2}}} \right).\]
Учитывая это получим
\[\mathop {\lim }\limits_{n \to  + \infty } \Delta \left( n \right) =  - \frac{1}{{{a^{\beta  - 1}}\Gamma \left( {2 - \beta } \right)}}.\]
Для завершения доказательства в качестве $\delta $  достаточно взять следующее число:
\[\delta  = \frac{1}{{{a^{\beta  - 1}}\Gamma \left( {2 - \beta } \right)}} - \varepsilon  > 0,\]
где  $0 < \varepsilon $ - произвольное малое число. \textbf{Лемма доказана.}

Перейдем к нахождению решения задачи $D$. Учитывая вышеизложенное, формальное решение задачи $D$ имеет вид
\[u\left( {x,y} \right) = \left\{ \begin{gathered}
  \sum\limits_{n = 1}^\infty  {{X_n}\left( x \right)\frac{{{\varphi _n}}}{{\Delta \left( n \right)}}{E_{\alpha ,1}}\left( { - {\lambda _n}{y^\alpha }} \right)} ,\,y > 0, \hfill \\
  \sum\limits_{n = 1}^\infty  {{X_n}\left( x \right)\left\{ {\frac{{{\varphi _n}}}{{\Delta \left( n \right)}}{E_{\beta ,1}}\left( { - {\lambda _n}{{\left( { - y} \right)}^\beta }} \right)} \right. + }  \hfill \\
  \left. { + \frac{{{\lambda _n}{\varphi _n}}}{{\Delta \left( n \right)}}\left( { - y} \right){E_{\beta ,2}}\left( { - {\lambda _n}{{\left( { - y} \right)}^\beta }} \right)} \right\},\,y < 0. \hfill \\
\end{gathered}  \right.\eqno(13)\]
Покажем , что (13) является классическим решением поставленной задачи. Справедлива теорема.

\textbf{Теорема 2.} Пусть  выполняются условия\\
1. $\Delta \left( n \right) \ne 0,\,\forall n;$\\
2. ${\varphi _1}\left( x \right) = {x^m}\varphi \left( x \right) \in {C^{2k}}\left( {{{\bar \Omega }_x}} \right),\varphi _1^{\left( j \right)}\left( 0 \right) = \varphi _1^{\left( j \right)}\left( 1 \right) = 0,$\\
3. ${\varphi _2}\left( x \right) = {\left( {{x^m}{\varphi _1}\left( x \right)} \right)^{\left( {2k} \right)}} \in {C^{2k}}\left( {{{\bar \Omega }_x}} \right),\varphi _2^{\left( j \right)}\left( 0 \right) = \varphi _2^{\left( j \right)}\left( 1 \right) = 0,$\\
4. ${\varphi _3}\left( x \right) = {\left( {{x^m}{\varphi _2}\left( x \right)} \right)^{\left( {2k} \right)}} \in {C^{2k}}\left( {{{\bar \Omega }_x}} \right),\varphi _3^{\left( j \right)}\left( 0 \right) = \varphi _3^{\left( j \right)}\left( 1 \right) = 0,$\\
5. ${x^{ - \frac{m}{2}}}\varphi _3^{\left( {2k} \right)}\left( x \right) \in {L_2}\left( {{{\overline \Omega  }_x}} \right),$
где $j = 0,1,...,k - 1.$\\
Тогда (13) является классическим решением задачи $D$.

\textbf{Доказательство. } Нужно показать возможность почленного дифференцирования ряда  (13) по переменным $x,y$  до порядков входящих в уравнение (1). Формально имеем
\[\frac{{{\partial ^{2k}}u\left( {x,y} \right)}}{{\partial {x^{2k}}}} = \left\{ \begin{gathered}
  {\left( { - 1} \right)^k}{x^{ - m}}\sum\limits_{n = 1}^\infty  {{\lambda _n}{X_n}\left( x \right)\frac{{{\varphi _n}}}{{\Delta \left( n \right)}}{E_{\alpha ,1}}\left( { - {\lambda _n}{y^\alpha }} \right)} ,\,y > 0, \hfill \\
  {\left( { - 1} \right)^k}{x^{ - m}}\sum\limits_{n = 1}^\infty  {{\lambda _n}{X_n}\left( x \right)\left\{ {\frac{{{\varphi _n}}}{{\Delta \left( n \right)}}{E_{\beta ,1}}\left( { - {\lambda _n}{{\left( { - y} \right)}^\beta }} \right)} \right. + }  \hfill \\
  \left. { + \frac{{{\lambda _n}{\varphi _n}}}{{\Delta \left( n \right)}}\left( { - y} \right){E_{\beta ,2}}\left( { - {\lambda _n}{{\left( { - y} \right)}^\beta }} \right)} \right\},\,y < 0. \hfill \\
\end{gathered}  \right.\]
Покажем сходимость ряда
\[\sum\limits_{n = 0}^\infty  {\frac{{\lambda _n^2{\varphi _n}{X_n}\left( x \right)}}{{\Delta \left( n \right)}}{E_{\beta ,2}}\left( { - {\lambda _n}{{\left( { - y} \right)}^\beta }} \right)} ,\]
сходимость остальных рядов показывается аналогично. В дальнейшем все независимые постоянные будем обозначать одной буквой $M$. Учитывая оценку [см.18]
\[\left| {{E_{\mu ,\eta }}\left( { - z} \right)} \right| \leqslant \frac{M}{{1 + \left| z \right|}},\,0 < \mu  < 2,M - const,\]
при $y<0$ имеем
\[\left| {{E_{\beta ,2}}\left( { - {\lambda _n}{{\left( { - y} \right)}^\beta }} \right)} \right| \leqslant M - const \Rightarrow \]
\[\left| {\sum\limits_{n = 1}^\infty  {\frac{{\lambda _n^2{\varphi _n}{X_n}\left( x \right)}}{{\Delta \left( n \right)}}{E_{\beta ,2}}\left( { - {\lambda _n}{{\left( { - y} \right)}^\beta }} \right)} } \right| \leqslant M\sum\limits_{n = 1}^\infty  {\lambda _n^2\left| {{X_n}\left( x \right)} \right|\left| {{\varphi _n}} \right|}  \leqslant \]
\[ \leqslant M\sum\limits_{n = 1}^\infty  {\frac{{\left| {{X_n}\left( x \right)} \right|}}{{{\lambda _n}}}\lambda _n^3\left| {{\varphi _n}} \right|}  \leqslant M\sqrt {\sum\limits_{n = 1}^\infty  {{{\left( {\frac{{\left| {{X_n}\left( x \right)} \right|}}{{{\lambda _n}}}} \right)}^2}} } \sqrt {\sum\limits_{n = 1}^\infty  {{{\left( {\lambda _n^3\left| {{\varphi _n}} \right|} \right)}^2}} }.\eqno(14)\]
Первый множитель сходится за счет оценки (9). Изучим второй множитель. Имеем
\[{\varphi _n} = \int\limits_0^1 {\varphi \left( x \right){X_n}\left( x \right)dx}  = \frac{{{{\left( { - 1} \right)}^k}}}{{{\lambda _n}}}\int\limits_0^1 {\varphi \left( x \right){x^m}X_n^{\left( {2k} \right)}\left( x \right)dx}  = \]
\[ = \frac{{{{\left( { - 1} \right)}^k}}}{{{\lambda _n}}}\left( {\left. {{\varphi _1}\left( x \right)X_n^{\left( {2k - 1} \right)}\left( x \right)} \right|_{x = 0}^{x = 1} - \int\limits_0^1 {{\varphi _1}\left( x \right)X_n^{\left( {2k - 1} \right)}\left( x \right)dx} } \right) = \]
\[ =  - \frac{{{{\left( { - 1} \right)}^k}}}{{{\lambda _n}}}\int\limits_0^1 {{\varphi _1}\left( x \right)X_n^{\left( {2k - 1} \right)}\left( x \right)dx}  = \frac{{{{\left( { - 1} \right)}^k}}}{{{\lambda _n}}}\int\limits_0^1 {\varphi _1^{\left( {2k} \right)}\left( x \right){X_n}\left( x \right)dx} .\]
Повторяя этот процесс два раза получим
\[{\varphi _n} = \frac{{{{\left( { - 1} \right)}^k}}}{{\lambda _n^3}}\int\limits_0^1 {\varphi _3^{\left( {2k} \right)}\left( x \right){X_n}\left( x \right)dx} .\]
Применим неравенство Бесселя
\[\sum\limits_{n = 1}^\infty  {{{\left( {\lambda _n^3{\varphi _n}} \right)}^2}}  \leqslant \int\limits_0^1 {{x^{ - m}}{{\left( {\varphi _3^{\left( {2k} \right)}\left( x \right)} \right)}^2}dx} .\eqno(15)\]
Из (9) и (15) следует  абсолютная и равномерная сходимость ряда (14). Итак выражение (13) является классическим решением задачи $D$. \textbf{Теорема 2 доказана}.

\textbf{Замечание .} Если ${\Delta }\left( n \right) = 0,$ при некоторых значениях $n = {n_1},{n_2},...{n_p}$, то для
разрешимости задачи $D$ достаточно выполнения условий ортогональности
 $\left( {{\varphi}\left( x \right),{X_n}\left( x \right)} \right) = 0,\,n = {n_1},...,{n_p}.$

\textbf{3. Единственость решения.} Справедлива теорема.

\textbf{Теорема 3.} Решение задачи $D$ единственно тогда и только тогда , когда
\[\Delta \left( n \right) = {E_{\alpha ,1}}\left( { - {\lambda _n}{b^\alpha }} \right) - \left( {{E_{\beta ,1}}\left( { - {\lambda _n}{a^\beta }} \right) + a{\lambda _n}{E_{\beta ,2}}\left( { - {\lambda _n}{a^\beta }} \right)} \right) \ne 0.\]

\textbf{Доказательство.} Пусть функция $u (x,y)$ есть решение задачи $D$ с однородными условиями. Рассмотрим его коэффициенты Фурье по системе собственных функций задачи (6)
\[{u_n}\left( y \right) = \int\limits_0^1 {{x^{ - m}}{X_n}\left( x \right)u\left( {x,y} \right)dx} = \]
\[ = \frac{{{{\left( { - 1} \right)}^k}}}{{{\lambda _n}}}\int\limits_0^1 {X_n^{\left( {2k} \right)}\left( x \right)u\left( {x,y} \right)dx}  = \]
\[ = \frac{{{{\left( { - 1} \right)}^k}}}{{{\lambda _n}}}\left( {\left. {X_n^{\left( {2k - 1} \right)}\left( x \right)u\left( {x,y} \right)} \right|_{x = 0}^{x = 1} - \int\limits_0^1 {X_n^{\left( {2k - 1} \right)}\left( x \right)\frac{{\partial u}}{{\partial x}}dx} } \right) = \]
\[ = \frac{{{{\left( { - 1} \right)}^k}}}{{{\lambda _n}}}\left( {\left. {{X_n}\left( x \right)\frac{{{\partial ^{2k - 1}}u}}{{\partial {x^{2k - 1}}}}} \right|_{x = 0}^{x = 1} - \int\limits_0^1 {{X_n}\left( x \right)\frac{{{\partial ^{2k}}u}}{{\partial {x^{2k}}}}dx} } \right) = \]
\[ = \frac{{{{\left( { - 1} \right)}^k}}}{{{\lambda _n}}}\left( {\left. {{X_n}\left( x \right)\frac{{{\partial ^{2k - 1}}u}}{{\partial {x^{2k - 1}}}}} \right|_{x = 0}^{x = 1} - \int\limits_0^1 {{X_n}\left( x \right)\frac{{{\partial ^{2k}}u}}{{\partial {x^{2k}}}}dx} } \right) = \]
\[ =  - \frac{1}{{{\lambda _n}}}\left\{ \begin{gathered}
  _CD_{0y}^\alpha {u_n}\left( y \right),\,y > 0, \hfill \\
  _CD_{0y}^\beta {u_n}\left( y \right),\,y < 0. \hfill \\
\end{gathered}  \right.\]
Итак для $u_{n}(y)$ получаем следующую задачу:
\[\left\{ \begin{gathered}
  _CD_{0y}^\alpha {u_n}\left( y \right) =  - {\lambda _n}{u_n}\left( y \right),{\mkern 1mu} {\kern 1pt} y > 0, \hfill \\
  _CD_{0y}^\beta {u_n}\left( y \right) =  - {\lambda _n}{u_n}\left( y \right),y < 0, \hfill \\
  {u_n}\left( { + 0} \right) = {u_n}\left( { - 0} \right), \hfill \\
  _CD_{0y}^\alpha {u_n}\left( { + 0} \right) = {{u'}_n}\left( { - 0} \right), \hfill \\
  {u_n}\left( b \right) = {u_n}\left( { - a} \right). \hfill \\
\end{gathered}  \right.\]
Теперь , если  $\Delta \left( n \right) = 0,$ при некотором $n$ , то задача (10) будет иметь ненулевое решение и отсюда свою очередь будет нарушаться единственность решения задачи $D$.\\
Пусть $\Delta \left( n \right) \ne 0,\,\forall n,$ тогда , как мы рассматривали выше , эта задача имеет только тривиальное решение
\[{u_n}\left( y \right) = 0,\,\forall n.\]
Далее имеем, т.к. ${\lambda _n} > 0,\forall n,\,\overline G \left( {x,\xi } \right) - $ симметрична, непрерывна и выполняются соотношения
$$\int\limits_0^1 {\overline G \left( {x,\xi } \right)dx}  < \infty ,\,\int\limits_0^1 {\overline G \left( {x,\xi } \right)d\xi }  < \infty ,\,\int\limits_0^1 {\int\limits_0^1 {\overline G \left( {x,\xi } \right)dx} d\xi }  < \infty ,$$
то из теоремы Мерсера следует, что
\[\overline G \left( {x,\xi } \right) = \sum\limits_{n = 1}^\infty  {\frac{{{{\overline X }_n}\left( x \right){{\overline X }_n}\left( \xi  \right)}}{{{\lambda _n}}}} .\]
Отсюда получим
\[{x^{ - \frac{m}{2}}}u\left( {x,y} \right) = \int\limits_0^1 {\overline G \left( {x,\xi } \right)\left( {{{\left( { - 1} \right)}^k}{\xi ^{\frac{m}{2}}}\frac{{{\partial ^{2k}}u\left( {\xi ,y} \right)}}{{\partial {\xi ^{2k}}}}} \right)d\xi }  = \]
\[ = {\left( { - 1} \right)^k}\int\limits_0^1 {\sum\limits_{n = 1}^\infty  {\frac{{{{\overline X }_n}\left( x \right){{\overline X }_n}\left( \xi  \right)}}{{{\lambda _n}}}} \left( {{\xi ^{\frac{m}{2}}}\frac{{{\partial ^{2k}}u\left( {\xi ,y} \right)}}{{\partial {\xi ^{2k}}}}} \right)} d\xi  = \]
\[ = {\left( { - 1} \right)^k}\sum\limits_{n = 1}^\infty  {\frac{{{x^{ - \frac{m}{2}}}{X_n}\left( x \right)}}{{{\lambda _n}}}} \int\limits_0^1 {\left( {{X_n}\left( \xi  \right)\frac{{{\partial ^{2k}}u\left( {\xi ,y} \right)}}{{\partial {\xi ^{2k}}}}} \right)} d\xi  = \]
т.к. ряд сходится равномерно, то можно поменять местами знаки интегрирования и суммы
\[ = {\left( { - 1} \right)^k}\sum\limits_{n = 1}^\infty  {\frac{{{x^{ - \frac{m}{2}}}{X_n}\left( x \right)}}{{{\lambda _n}}}} \int\limits_0^1 {X_n^{\left( {2k} \right)}\left( \xi  \right)u\left( {\xi ,y} \right)} d\xi  = \]
\[ = {x^{ - \frac{m}{2}}}\sum\limits_{n = 1}^\infty  {{X_n}\left( x \right)} \int\limits_0^1 {{\xi ^{ - \frac{m}{2}}}{X_n}\left( \xi  \right)u\left( {\xi ,y} \right)} d\xi  = 0.\]
\textbf{Теорема 3} доказана.

\begin{center}
Литература
\end{center}
1.  А. М. Нахушев . Дробное исчисление и его применение // М., Физматлит,2003, с. 272.\\
2.  А. Н. Боголюбов , А. А. Кобликов , Д. Д. Смирнова , Н. Е. Шапкина.  Математическое моделирование сред с временной дисперсией при помощи дробного дифференцирования // Матем. моделирование, 25: 12,  2013, с. 50–64.\\
3. И.М. Гельфанд.  Некоторые вопросы анализа и дифференциальных ураванений // УМН. 1959. Т. XIV. вып. 3(87). С. 3-19.\\
4.	Г.М. Стручина.  Задача о сопряжении двух уравнений // Инженерно-физический журнал. 1961. Т. IV. № 11. С. 99-104.\\
5.	 Я.С. Уфлянд.  К вопросу о распространении колебаний в составных электрических линиях // Инженерно-физический журнал. 1964. Т. VII. № 1. С. 89-92.\\
6. К. Б. Сабитов, Начально-граничная и обратные задачи для неоднородного уравнения смешанного параболо-гиперболического уравнения // Матем. заметки, 2017,том 102, выпуск 3, 415–435\\
7. К. Б. Сабитов, Нелокальная задача для уравнения параболо-гиперболического типа в прямоугольной области // Матем. заметки, 2011, том 89, выпуск 4, 596–602.\\
8. К. Б. Сабитов, Э. М. Сафин, Обратная задача для уравнения смешанного параболо-гиперболического типа // Матем. заметки, 2010, том 87, выпуск 6, 907–918.\\
9. К. Б. Сабитов, Задача Трикоми для уравнения смешанного параболо-гиперболического в прямоугольной области // Матем. заметки, 2009, том 86, выпуск 2, 273–279\\
10. A.S. Berdyshev, A. Cabada, E.T. Karimov. On a non-local boundary problem for a parabolic–hyperbolic equation involving a Riemann–Liouville fractional differential operator // Nonlinear Analysis: Theory, Methods and Applications. 75 (6), 2012, p. 3268-3273.\\
11. P.Agarwal, A.S. Berdyshev, E. T. Karimov. Solvability of a non-local problem with integral transmitting condition for mixed type equation with Caputo fractional derivative // Results in Mathematics. 71(3), 2017, p.1235-1257.\\
12. A.S. Berdyshev, E.T Karimov, N. Akhtaeva. Boundary value problems with integral gluing conditions for fractional-order mixed-type equation // International Journal of Differential Equations. 2011.\\
13. O. Kh. Abdullaev , K. Sadarangani  Non-local problems with integral gluing condition for loaded mixed
type equations involving the Caputo fractional derivative // Electron. J. Differential Equations, 2016.\\
14. O. Kh. Masaeva. Uniqueness of solutions to Dirichlet problems for generalized Lavrent’ev-Bitsadze
equations with a fractional derivative // Electron. J. Differential Eq., 2017 (74): 1–8.\\
15. A. S. Berdyshev , A. Cabada , B. J. Kadirkulov.  The Samarskii–Ionkin type problem for the fourth order
parabolic equation with fractional differential operator // Comput. Math. Appl., 2011. Vol. 62. P. 3884–
3893. DOI: 10.1016/j.camwa.2011.09.038 Vol. 2016. No. 164. P. 1–10. URL: https://ejde.math.txstate.edu\\
16. A. S. Berdyshev , B. J. Kadirkulov . On a nonlocal problem for a fourth-order parabolic equation
with the fractional Dzhrbashyan–Nersesyan operator // Differ. Equ., 2016. Vol. 52. No. 1. P. 122–127.
DOI: 10.1134/S0012266116010109\\
17. Yu. P. Apakov, B.Yu. Irgashev. Boundary-value problem for a degenerate high-odd-order equation. Ukrainian Mathematical Journal.66(10)(2015), p. 1475-1490.\\
18. М. М. Джрбашян .  Интегральные преобразования и представления функций в комплексной
области // M.: Наука.1966, с. 672.\\
19. А. В. Псху.  О вещественных нулях функции типа Миттаг – Леффлера // Матем. заметки, 77 (4),2005, с.592–599.\\

\end{document}